\documentclass[12pt]{article}

\def\Dj{\hbox{D\kern-.73em\raise.30ex\hbox{-}
\raise-.30ex\hbox{}}}
\def\dj{\hbox{d\kern-.33em\raise.80ex\hbox{-}
\raise-.80ex\hbox{\kern-.40em}}}

\usepackage{amssymb}
\usepackage{amsmath} % If you want to use AMS-LaTeX, comment out these
\usepackage{cite}
\usepackage[dvips]{graphicx}
\usepackage[dvipdfm]{color}
\usepackage[dvipdfm,colorlinks=true,citecolor=blue,linkcolor=red]{hyperref}

\newenvironment{Proof}
        {\par\noindent{\bf Proof.}}
        {\hspace*{1mm}{$\Box$}\\[1mm]}
\begin{document}
\newtheorem{Theorem}{Theorem}[section]
\newtheorem{Thm}[Theorem]{Theorem}
\newtheorem{Lem}[Theorem]{Lemma}
\newtheorem{Con}[Theorem]{Conjecture}
\newtheorem{Cor}[Theorem]{Corollary}
\newtheorem{Pro}[Theorem]{Proposition}
\newtheorem{Conj}[Theorem]{Conjecture}
\newtheorem{Def}[Theorem]{Definition}
\newtheorem{Que}[Theorem]{Question}
\newtheorem{ex}[Theorem]{Example}
\newcommand{\eq}{\begin{equation}}
\newcommand{\en}{\end{equation}}

\def\G{\overrightarrow{G}}
\def\a{\lambda}
\def\g1{{\G}^{+}}
\def\g2{{\G}^{-}}
\def\zm{\noindent{\bf Proof.\ }}

\begin{center}

\thispagestyle{empty}
{\Large Oriented unicyclic  graphs with extremal skew energy}
\vspace*{8mm}\\

     {Yaoping Hou,\ Xiaoling Shen and \ Chongyan Zhang }
   \\[2mm]
   { Department of Mathematics, Hunan Normal
             University\\ Changsha,  Hunan 410081, China\\ email: {\tt
yphou@hunnu.edu.cn}}\\[6mm]

%(Received  May 2011 )
\end{center}

\baselineskip=0.20in

{\bf Abstract}

Let $\G$ be an oriented graph of order $n$ and $\a_1,\a_2,\cdots, \a_n$ denote all the eigenvalues of the skew-adjacency matrix of $\G.$ The skew energy  $\displaystyle{\cal E}_s(\G)= \sum_{i=1}^{n} |\a_i|.$
  In this paper,
the  oriented unicyclic  graphs with minimal and maximal skew energy are
determined.


AMS Classification  05C50, 15A03

 \baselineskip=0.25in

 \section{Introduction }

\ \ \ \ \ \ An important quantum--chemical characteristic of a conjugated
molecule is its total $\pi$-electron energy.
The energy of a graph  has closed links to chemistry. Since the concept of the energy of simple undirected
graphs was introduced by Gutman in \cite{gutman}, there have been  lots
of research papers   on this  topic.  For the extremal energy of  unicyclic graph, Hou
\cite{hou1} showed that $S^3_ n$ is the graph with minimal energy
in all unicyclic graphs;  In \cite{bh}, Huo and Li  showed that
$P^6_n$ is the graph with maximal energy in all
unicyclic graphs. More results on  the energy of unicyclic graphs see
\cite{gu,hh2,hou2,lx,whw,whw2}.

There are various  generalizations of energy to   graph matrices \cite{vn}, such as Laplacian energy, incidence energy and distance energy.
 Let $G$ be a simple undirected finite graph of  order
$n$  and  $\overrightarrow{G}$ be
an orientation of $G,$
 which assigns to each edge of $G$ a direction  so that $\G$ becomes a directed graph. {\em All digraphs in this paper are oriented graphs of some graphs.}  Let $\G$ be an oriented graph of order $n$.  {\em The skew adjacency matrix $S(\G)=(s_{i,j})$ }is a real skew symmetric matrix,  where $s_{i,j}=1 $ and $  s_{j,i}=-1$ if  $i\rightarrow j$ is an arc of $\G,$ otherwise $s_{i,j}=s_{j,i}=0.$ The {\em skew spectrum $Sp(\G)$} of $\G$ is defined as the spectrum of $S(\G) .$  Note that $Sp(\G) $ consists of only purely  imaginary eigenvalues because $S(\G)$ is real skew symmetric.

Recently, the skew-energy of  an  oriented graph $\G$   was defined as  the energy of matrix $S(\G)$ in \cite{ca}, that is,
$${\cal E}_s(\G)= \sum_{\a \in Sp(\G)} |\a|.$$

There are situations when chemists use digraphs rather than graphs. One such situation is when vertices represent distinct chemical species and arcs represent the direction in which a particular reaction takes place between the two corresponding species. It is possible that the skew energy has similar applications  as energy in chemistry.   For a graph  $G$, there are any orientations on it,   it is also interesting to find what orientation  has extremal energy among all orientations of a given graph.

 An unicyclic graph is  the connected graph with the same number of vertices and edges.  
 In this paper,
we are interested in studying the  orientations of  unicyclic  graphs with extremal
skew energy. Let $G(n,\ell)$ be the set of all connected unicyclic
graphs on $n$ vertices with girth $\ell$. Denote, as usual, the
$n$-vertex path and cycle by $P_n$ and $C_n$, respectively. Let
$P^{\ell}_n$ be the unicyclic graph obtained by connecting a vertex of
$C_{\ell}$ with a terminal vertex of $P_{n-\ell}$, $S_n^{\ell}$ be the graph
obtained by connecting $n-\ell$ pendant vertices to a vertex of $C_{\ell}$
(see Fig. $1$).

The rest of the paper is organized as follows.  In section $2$, a
new integral formal for ${\cal E}_s(\overrightarrow{G})$ is
obtained and the oriented  unicyclic graph  with minimal skew energy is
determined. In section $3$, the  oriented unicyclic graph with maximal skew energy
is determined.

\section{Oriented unicyclic graph with minimal skew energy}

\ \ \ \ \ \ \ Let $G$ be a graph. A \emph{linear subgraph} $L$ of
$G$ is a disjoint union of some edges and some cycles in $G$
(\cite{dmc}). A \emph{$k$-matching} $M$ in $G$ is a disjoint union
of $k$-edges. If $2k$ is the order of $G$, then a $k$-matching of
$G$ is called a \emph{perfect matching }of $G$. The number of
$k$-matching is denoted by $m({G},k)$.

If $C$ be any undirected even cycle of ${G},$  we
say $C$ is\emph{ evenly oriented } relative to  the orientation
$\overrightarrow{G}$ if it has an even number of edges oriented in
 clockwise direction. Otherwise $C$ is \emph{oddly
oriented}.

We call a linear subgraph $L$ of $G$  \emph{evenly linear} if $L$
contains no  odd cycle  and denote by $\mathcal{E}
\mathcal{L}_i(G)$ (or $\mathcal{E} \mathcal{L}_i$ for short)
the set of all evenly linear subgraphs of $G$ with $i$ vertices.
For a linear subgraph $L\in \mathcal{E} \mathcal{L}_i$ denote by
$p_e(L) $(resp., $p_o(L)$ )  the number of evenly (resp., oddly)
oriented cycles in $L$ relative to $\overrightarrow{G}$. Denote
the characteristic polynomial of $S(\overrightarrow{ G} )$ by

\eq P_s(\overrightarrow{ G }; x)=det(xI -S(\overrightarrow{ G }))
= \sum ^{ n }_{i=0} b_ix^{n-i}.\label{ps}\en

Then (i) $b_0 = 1$, (ii) $b_2$ is the number of edges of $G,$ (iii) all $b_i\geq 0$
and (iv) $b_i=0$ for all odd $i$
since the determinant of every real skew  symmetric matrix is nonnegative and is 0 if its order is odd.

\includegraphics[width=13.8cm]{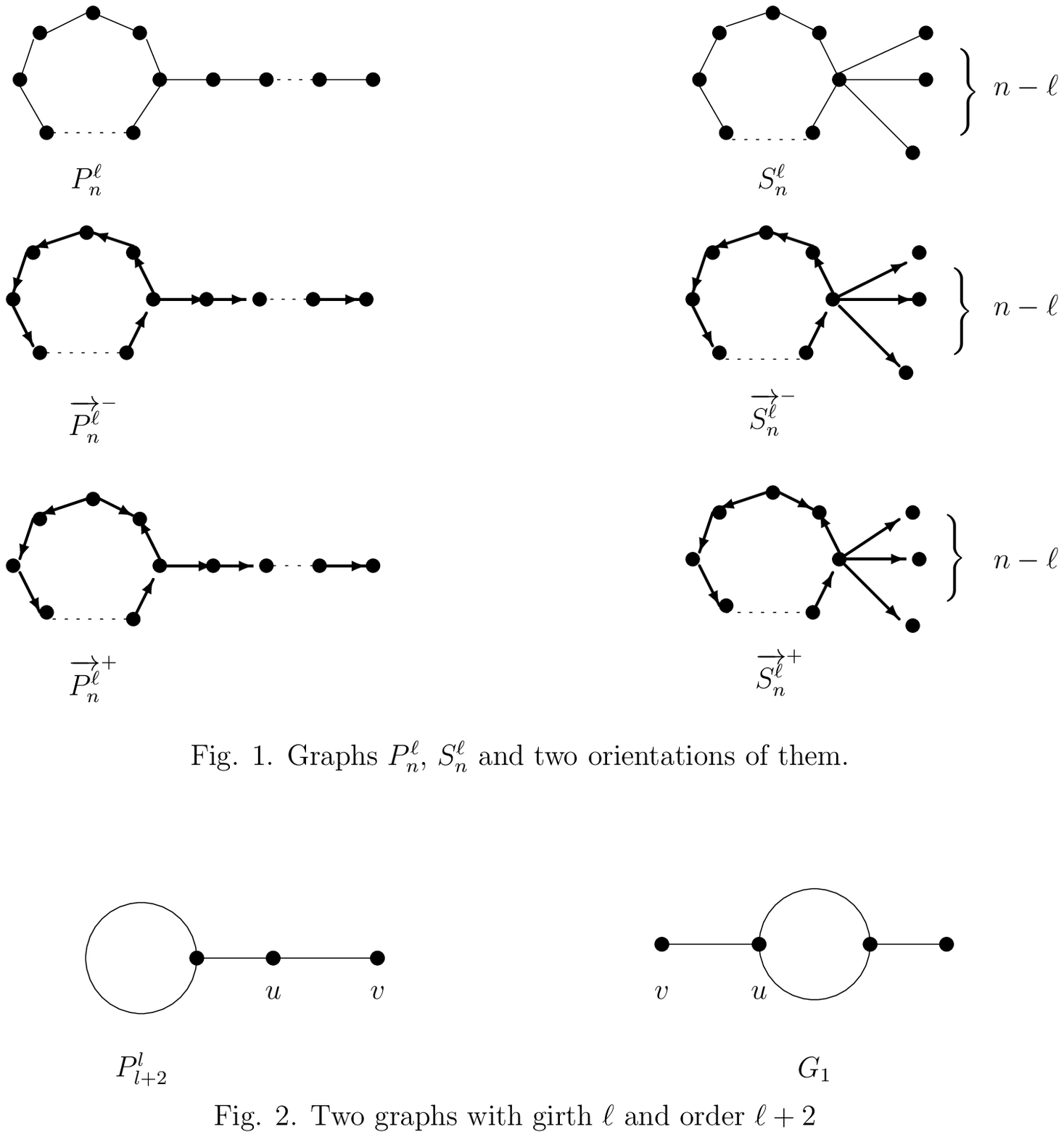}

\vspace*{-5cm}

We have

\begin{Lem} (\cite{hou}) Let
$\overrightarrow{G}$ be an orientation of  a graph $G$. Then \eq
b_i(\overrightarrow{G})=\sum _{L\in \mathcal{E} \mathcal{L}_i}
(-2)^{p_e(L)}2^{p_o(L)},\label{b2k}\en
 where $p_e(L)$ is the number of evenly oriented cycles
of $L$  and $p_o(S)$ is the number
of oddly oriented cycles of $L$ relative to $\overrightarrow{G}$,
respectively.\label{char}
\end{Lem}

From Lemma \ref{char}, we can get the following statement.

\begin{Cor}  Let  $G$ be a unicyclic graph with  unique cycle $C$ of length $\ell$ and $\G$ be an orientation of $G.$  Then
$$b_{2k}(\overrightarrow{G})=\left\{\begin{array}{cc}
m({G},k),&
   C\ is\  an\ odd\  cycle;\\
m({G},k)-2m({G}-C,k-\frac{\ell}{2}),& C\ is \  evenly\  oriented;\\
m({G},k)+2m({G}-C,k-\frac{\ell}{2}),& C\
\ is\ oddly\  oriented.
\end{array}\right. $$  \label{coef2}\end{Cor}

\begin{Lem} (\cite{hou}) Let $e=uv$ be an edge of $G$ that is on no any even cycle of $G.$ Then
$$P_s(
\overrightarrow{G}; x) = P_s(\overrightarrow{G}-e; x) +
P_s(\overrightarrow{G}-u-v; x).
$$\label{cut}
\end{Lem}

By Lemma \ref{cut}, we obtain that

\begin{Cor}Let $e=uv$   be an edge of $G$ that is on no  even cycle of $G.$
Then \eq
b_{2k}(\overrightarrow{G})=b_{2k}(\overrightarrow{G}-e)+b_{2k-2}(\overrightarrow{G}-u-v).\label{brela}\en

 Furthermore, if $e=uv$ is a pendant  edge with the pendant vertex
$v.$ Then
 \eq
b_{2k}(\overrightarrow{G})=b_{2k}(\overrightarrow{G}-v)+b_{2k-2}(\overrightarrow{G}-u-v).\label{brela2}\en
\label{bre}\end{Cor}

For  any orientation of  a graph which does not contain  any even cycle (in particular, a tree,  a unicyclic non-bipartite graph) by Theorem \ref{coeff}, we have
$b_{2k}(\overrightarrow{G})=m({G},k)$.

 For the $k$-matching number of a graph $G$,   the following result is well-known.

\begin{Lem} Let $e=uv$ be an edge of $G$. Then\\
(i) $
m({G},k)=m({G}-e,k)+m({G}-u-v,k-1).$\\
(ii) If $G$ is a forest, then
$m({G},k)\leq m(P_n,k)$, $k\geq 1$.\\
(iii) If ${H}$ is a subgraph of
${G}$, then $m({H},k)\leq
m({G},k)$, $k\geq 1$. Moreover, if
${H}$ is a proper subgraph of ${G}$,
then the inequality is strict. \label{match}\end{Lem}

 From
\cite{ca}, an integral formula for skew energy was given: \eq
{\cal E}_s(\overrightarrow{G})= \frac{1}{\pi}
\int^{\infty}_{-\infty}( n + x
\frac{P'_s(\overrightarrow{G},-x)}{P_s(\overrightarrow{G},-x)})dx,\label{1}\en
 where $P'_s(\overrightarrow{G},-x)$ is the derivative of $P_s(\overrightarrow{G},-x).$

\begin{Thm}  Let $\overrightarrow{G}$ be an orientation of  a graph $G$. Then
 $${\cal E}_s(\overrightarrow{G})= \frac{1}{\pi} \int^{\infty}_{-\infty}\frac{1}{t^2}ln(1+\sum_{k=1}^{\lfloor
 \frac{n}{2}\rfloor}
 b_{2k}t^{2k})dt.$$\label{formu} \end{Thm}

{\bf Proof.} From equality (\ref{1}), we have
\begin{eqnarray*}
 {\cal E}_s(\overrightarrow{G})&=& \frac{1}{\pi}
\int^{\infty}_{-\infty}( n + x
\frac{P'_s(\overrightarrow{G},-x)}{P_s(\overrightarrow{G},-x)})dx
= \frac{1}{\pi} \int^{\infty}_{-\infty}( n
-\mu\frac{P'_s(\mu)}{P_s(\mu)})d\mu\\
&=&\frac{1}{\pi} \int^{\infty}_0( n
-\frac{1}{t}\frac{P'_s(\frac{1}{t})}{P_s(\frac{1}{t})})(-\frac{1}{t^2})dt+
\frac{1}{\pi} \int^0_{-\infty}( n
-\frac{1}{t}\frac{P'_s(\frac{1}{t})}{P_s(\frac{1}{t})})(-\frac{1}{t^2})dt\\
&=&\frac{1}{\pi} \int^{\infty}_{-\infty}( n
-\frac{1}{t}\frac{P'_s(\frac{1}{t})}{P_s(\frac{1}{t})})\frac{1}{t^2}dt
=\frac{1}{\pi} \int^{\infty}_{-\infty} \frac{1}{t}d[ln(t^nP_s(\frac{1}{t})]\\
&=&\frac{1}{\pi} \int^{\infty}_{-\infty}
\frac{1}{t^2}ln(t^nP_s(\frac{1}{t}))dt.
= \frac{1}{\pi} \int^{\infty}_{-\infty}\frac{1}{t^2}ln(1+\sum_{k=1}^{\lfloor
 \frac{n}{2}\rfloor}
 b_{2k}t^{2k})dt.  \;\; \; \Box
\label{energy}
\end{eqnarray*}

\begin{ex} {\rm Let $\overrightarrow{C_4}^{-}$ be an orientation of $C_4$ such that all  edges have the same direction.
Then by Lemma \ref{char}, the characteristic polynomial of
$S(\overrightarrow{C_4}^{-})$ is  $x^4+4x^2$. By Lemma \ref{formu},
\begin{eqnarray*}{\cal E}_s(\overrightarrow{C_4}^{-},x))&=&\frac{1}{\pi} \int^{\infty}_{-\infty}\frac{1}{t^2}ln(1+4t^2)dt\\
&=&
\frac{8}{\pi}\int^{\infty}_{0}\frac{1}{1+4t^2}dt=4.
\end{eqnarray*}}
\end{ex}

From Theorem \ref{formu}, ${\cal E}_s(\overrightarrow{G})$ is
an
 increasing function of $b_{2k}(\overrightarrow{G})$, $k = 0$, $1$, $\cdots$, $\lfloor\frac{n}{2}\rfloor$.
Consequently, if $\overrightarrow{G_1}$ and $\overrightarrow{G_2}$
are  oriented graphs of $G_1$ and $G_2,$ respectively, for which \eq
b_{2k}(\overrightarrow{G_1})\geq
b_{2k}(\overrightarrow{G_2})\label{b}\en for all
$\lfloor{\frac{n}{2}}\rfloor\geq k\geq 0$, then

\eq {\cal E}_s(\overrightarrow{G_1})\geq
{\cal E}_s(\overrightarrow{G_2})\label{e}.\en

Equality in (\ref{e}) is attained only if (\ref{b}) is an
equality for all  $\lfloor{\frac{n}{2}}\rfloor\geq k\geq 0$.

If relations (\ref{b}) hold for all $k$, then we write $G_1\succeq
G_2$ or $G_2\preceq G_1$. If $G_1\succeq G_2$, but not $G_2\succeq
G_1$, then we write $G_1\succ G_2$.

Let $\overrightarrow{G}$ be an orientation of a graph $G$. Let $W$
be a subset of $V(G)$ and $\overline{W}=V(G)\setminus W$. The
orientation
 $\overrightarrow{G'}$ of $G $ obtained from $\overrightarrow{G}$
  by reversing the direction of
all arcs between $\overline{W}$ and $W$ is said to be obtained
from $\overrightarrow{G}$ by switching with respect to $W$.
Moreover, two orientations $\overrightarrow{G}$ and
$\overrightarrow{G'}$ of a graph $G$ are said to be
switching-equivalent if $\overrightarrow{G'}$ can be obtained from
$\overrightarrow{G}$ by a sequences of switching.  If two orientations $\G$ and $\G '$ of a graph  $G$ are switching-equivalent then their skew-adjacency matrices  are  similar  by a diagonal matrix whose $(i,i)$-entry is $-1$  when $i\in W $ and $1 $ when $i \not \in W$  by Lemma 3.1 in \cite{ca}  and hence  they have the same  skew spectrum. Thus,

\begin{Lem}  Let $\G$ and $\G'$ be two orientations of a graph $G.$  If $\G$ and $\G'$ are switching-equivalent,   then  ${\cal E}_s(\G) =
{\cal E}_s(\G')$.\label{switch}\end{Lem}

By Lemma \ref{switch} and switching-equivalence, there are only two different
orientations on a unicyclic graph $G.$ All  edges on the
unique cycle $C$ have the same direction or just one edge on the
cycle has the opposite direction to the directions of other
edges on the cycle regardless how the edges not on the cycle $C$
are oriented. Denote by ${\G}^-$ (${\G}^+,$ resp.) the orientation
of $G$ in first (second,  resp.) case above.

\begin{Lem} For any unicyclic graph $G$, ${\G}^+\succeq {\G}^-.$ \label{coeff}    \end{Lem}

\begin{Proof} By Corollary  (\ref{coef2}), if the girth $\ell$ of $G$ is odd, then $$b_{2k}({{\G}^+})=
b_{2k}({{\G}^-})=m(G,k)  \mbox{ for all $0\leq k \leq
\lfloor\frac{n}{2}\rfloor$.}$$ If the girth $\ell$ of $G$  is even, then
$$b_{2k}({{\G}^+})=m(G, k)+
2m(G-C_{\ell},k-\lfloor\frac{\ell}{2}\rfloor),$$
$$b_{2k}({{\G}^-})=m(G, k)-
2m({G}-C_{\ell},k-\lfloor\frac{\ell}{2}\rfloor)$$ for all
nonnegative integer $\lfloor\frac{n}{2}\rfloor \geq k \geq
\lfloor\frac{\ell}{2}\rfloor,$  and
$$b_0({{\G}^+})=b_0({{\G}^-})=1, \;\;
b_{2k}({{\G}^+})=b_{2k}({{\G}^-})=m(G,k)$$ for all $0<
k<\lfloor\frac{\ell}{2}\rfloor$. Thus the result follows immediately.
\end{Proof}

\begin{Lem} Let ${\G}$ be an orientation of a unicyclic graph $G\in G(n,\ell)$,
$G\neq S_n^{\ell}$.  If  unique cycle $C_{\ell}$  in ${\G}$ and $\overrightarrow{S_{n}^{\ell}}$  is the same orientation, then
${\G}\succ \overrightarrow{S_{n}^{\ell}}$.\label{sn1}
\end{Lem}

\begin{Proof} We prove the statement by induction on $n$. Since $G\neq S_n^{\ell}$, $n\geq \ell+2 $.  For $n=\ell+2$,
$G$ is one of the two graphs in Fig. $2$.

\noindent By Corollary  \ref{bre},
$$b_{2k}(\overrightarrow{P_{\ell+2}^{\ell}})=b_{2k}(\overrightarrow{P_{\ell+1}^{\ell}})
+b_{2k-2}(\overrightarrow{C_{\ell}}),$$

$$b_{2k}({\G_1})=b_{2k}(\overrightarrow{P_{\ell+1}^{\ell}})+
b_{2k-2}(\overrightarrow{T}),\label{t}$$

\noindent where $T$ is a graph obtained by connecting one of the
vertices in $P_{\ell-1}$ to a pendant vertex.

$$b_{2k}(\overrightarrow{S_{\ell+2}^{\ell}})=b_{2k}(\overrightarrow{P_{\ell+1}^{\ell}})+
b_{2k-2}(\overrightarrow{P_{\ell-1}}).$$

\noindent By Corollary \ref{bre}, if $\ell$ is  odd or $\ell$ is even
but $k\leq \frac{\ell}{2}$  then
$$b_{2k-2}(\overrightarrow{C_{\ell}})=m({C_{\ell}},k-1);$$ for
$k= \frac{\ell}{2}+1.$  If $C$ is oddly oriented, then
$$b_{l}(\overrightarrow{C_{\ell}})=m({C_{\ell}},\frac{\ell}{2})+2=4;$$
If $C$ is evenly oriented, then
$$b_{\ell}(\overrightarrow{C_{\ell}})=m({C_{\ell}},\frac{\ell}{2})-2=
0,\;\; b_{\ell}(\overrightarrow{P_{\ell-1}})=0.$$ Since both $T$ and
$P_{\ell-1}$ are trees,
$$b_{2k-2}(\overrightarrow{T})=m({T},k-1),$$
$$b_{2k-2}(\overrightarrow{P_{\ell-1}})=m({P_{\ell-1}},k-1).$$
Since $P_{\ell-1}$ is a proper subgraph of both  $C_{\ell}$ and $T$, by Lemma \ref{match},
$$m(\overrightarrow{P_{\ell-1}},k-1)< m({T},k-1),
m({P_{\ell-1}},k-1)< m({C_{\ell}},k-1).$$
The result holds immediately for $n=\ell+2$.

 Suppose that
${\G}\succ \overrightarrow{S_{n'}^{\ell}}$ for all $n'<
n$. Since ${\G}$ is a unicyclic digraph, there is at
least a pendant edge $uv$  with pendant vertex $v$ in
${\G}$, by equality (\ref{brela2}), we get

\begin{eqnarray*}
b_{2k}({\G})&=&b_{2k}({\G}-v)+
b_{2k-2}({\G}-u-v), \\
b_{2k}(\overrightarrow{S_n^{\ell}})&=& b_{2k}
(\overrightarrow{S_{n-1}^{\ell}})
+b_{2k-2}(\overrightarrow{P_{\ell-1}}).\end{eqnarray*}

By induction assumption, it  suffices to prove that
$b_{2k-2}({\G}-u-v)\geq
b_{2k-2}(\overrightarrow{P_{\ell-1}})$ for all $ 0\leq k\leq
\lfloor\frac{n-\ell}{2}\rfloor$. For $k>
\lfloor\frac{\ell-1}{2}\rfloor$, we have
$b_{2k-2}({\G}-u-v)\geq
b_{2k-2}(\overrightarrow{P_{\ell-1}})=0$. For
$\lfloor\frac{\ell-1}{2}\rfloor \geq k\geq 0 $,
$b_{2k-2}({\G}-u-v)=
m({\G}-u-v,k-1)\geq
m(\overrightarrow{P_{\ell-1}},k-1)=b_{2k-2}(\overrightarrow{P_{\ell-1}})$
since ${P_{\ell-1}}$ is a subgraph of
${G}-u-v$.
\end{Proof}

\begin{Lem} Let  $n \geq \ell\geq 6$ or $n> \ell=5$, then
${\overrightarrow{S_n^4}}^-\prec{\overrightarrow{S_n^4}}^+
\prec{\overrightarrow{S_n^{\ell}}}^-
\preceq{\overrightarrow{S_n^{\ell}}}^+$.\label{s4}
\end{Lem}

\begin{Proof}  By Lemma
\ref{char}, the characteristic polynomial of
$\overrightarrow{S_n}^4$ is:

\begin{eqnarray*}
P_s({\overrightarrow{S_n^4}}^{+},x)&=&x^{n-4}(x^4+nx^2+2n-4);\\
P_s({\overrightarrow{S_n^{4}}}^-,x)&=&x^{n-4}(x^4+nx^2+2n-8);
\end{eqnarray*}
 Obviously, ${\overrightarrow{S_n^4}}^-\prec{\overrightarrow{S_n^4}}^+$.

So it  suffices to prove that $b_4(\overrightarrow{S_n^{\ell}})>
2n-4$. By equality (\ref{b2k}), \begin{eqnarray*}
b_4({\overrightarrow{S_n^{\ell}}}^+)&= &b_4({\overrightarrow{S_n^{\ell}}}^-)\\
&=&m(P_{\ell-1},2)+(n-\ell)m(P_{\ell-1},1)+2m(P_{\ell-2},1)\\
&=&
\frac{(\ell-3)(\ell-4)}{2}+(n-\ell)(\ell-2)+2(\ell-3)\\
&=&\frac{2n\ell-\ell
^2+\ell-4n}{2}.
\end{eqnarray*}

For $n>\ell=5$, $b_4(\overrightarrow{S_n^{\ell}})=3n-10\geq 2n-4$; For
$n\geq \ell=6$, $b_4(\overrightarrow{S_n^{\ell}})=4n-15>2n-4$; for $n\geq
\ell=7$, $b_4(\overrightarrow{S_n^{\ell}})=5n-21>2n-4$; for $\ell\geq 8$, $
b_4(\overrightarrow{S_n^{\ell}})\geq \frac{n\ell+\ell-4n}{2}\geq 2n+4>2n-4.$

By Lemma \ref{coeff}, the proof is completed.
\end{Proof}

For $n=\ell=5$,
$b_4({\overrightarrow{S_n^{\ell}}}^+)=\frac{2n{\ell}-{\ell}^2+\ell-4n}{2}=5$,
$b_4({\overrightarrow{S_n^{4}}}^+)=10-4=6$,
$b_4({\overrightarrow{S_n^{4}}}^-)=10-8=2$. It gives
${\overrightarrow{S_5^{4}}}^-\prec
{\overrightarrow{S_5^{5}}}^+\prec {\overrightarrow{S_5^{4}}}^+$.

\begin{Lem} Let $\overrightarrow{S_n^3}$ be any orientation of unicyclic graph $S_n^3$.  Then
 $\overrightarrow{S_n^3}\prec \overrightarrow{S_n^4}$  for $n\geq 6$;  $\overrightarrow{S_n^3}={\overrightarrow{S_n^{4}}}^-\prec {\overrightarrow{S_n^4}}^+$ for $n=5$;
  ${\overrightarrow{S_n^{4}}}^-\prec \overrightarrow{S_n^{3}}\prec {\overrightarrow{S_n^{4}}}^+$ for $n=4$.\label{s3}\end{Lem}

\begin{Proof} By Lemma \ref{char}, the characteristic polynomial of $\overrightarrow{S_n^3}$ is
$$P_s(\overrightarrow{S_n^3},x)=x^{n-4}(x^4+nx^2+n-3).$$ Since $n-3<
2n-8$, the result holds by Lemmas \ref{s4} and \ref{coeff}.

For $n=5$, it is easy to get that for any orientation ${\cal
E}_s(\overrightarrow{S_5^3})=2\sqrt{5+2\sqrt{2}}$, ${\cal
E}_s({\overrightarrow{S_5^{4}}}^-)=2\sqrt{5+2\sqrt{2}}<
2(\sqrt{2}+\sqrt{3})$, ${\cal
E}_s({\overrightarrow{S_5^{4}}}^+)=2(\sqrt{2}+\sqrt{3}) $. Then
$\overrightarrow{S_n^3}={\overrightarrow{S_5^{4}}}^-\prec
{\overrightarrow{S_5^{4}}}^+$.

 For $n=4$,  ${\cal E}_s(\overrightarrow{S_5^3})=2\sqrt{6}$,
 ${\cal E}_s({\overrightarrow{C_4}}^+)=4\sqrt{2}$,
 ${\cal E}_s({\overrightarrow{C_4}}^-)=4$.
  Thus ${\overrightarrow{S_5^{4}}}^-\prec
\overrightarrow{S_n^3}\prec
{\overrightarrow{S_5^{4}}}^+$.\end{Proof}

From Lemmas \ref{sn1}, \ref{s4} and  \ref{s3},
 we obtain one of the main result of this paper.

 \begin{Thm} Among all orientations of unicyclic graphs on $n$ vertices,  $ \overrightarrow{S_n^3}$ has the minimal skew energy   and
 ${\overrightarrow{S_n^4}}^-$ has the second minimal skew energy for $n\geq 6$; both $ {\overrightarrow{S_5^3}}$ and $ {\overrightarrow{S_5^4}} ^-$
 have the minimal skew energy, ${\overrightarrow{S_5^4}} ^+$ has the second minimal skew energy  for $n=5$;   ${\overrightarrow{C^4}}^-$ has the minimal skew energy,
  $\overrightarrow{S_4^3}$ has the second  minimal skew energy for $n=4$. \end{Thm}

\section{ Oriented unicyclic graph with maximal skew energy}

 By Lemma
\ref{coeff}, we only need to consider ${{\G}^+}$ for considering of maximum skew energy.

\begin{Lem} Let ${\G}$ be an orientation of unicyclic graph $G\in G(n,\ell)$ and  $G\neq P_n^{\ell}$. Then ${{\G}^+}\prec
{\overrightarrow{P_n^{\ell}}}^+$. \label{gless}\end{Lem}

 \zm  We prove the statement by induction on $n$. For $n=\ell+2$,  there are only two cases for
$G\neq P_n^l$: one is $S_{\ell+2}^l$, the other is the graph $G_1$ in
Fig. $2$. By the proof in  Lemma \ref{sn1}, we only need to prove
that $b_{2k-2}(\overrightarrow{C_l})\geq
m({T},k-1)$. By Lemma \ref{match},
\begin{eqnarray*}
m({\overrightarrow{C_{\ell}}}^+,k-1)&= & m({C_{\ell}},k-1)\\
&=& m({P_{\ell}},k-1)+m({P_{\ell-2}},k-2),\\
 m({T},k-1)&=&m({P_{\ell-1}},k-1)+
m({P_{s}}\bigcup {P_t},k-2)\\
&\leq &
m({P_{\ell-1}},k-1)+m({P_{\ell-2}},k-2)\\
&<&
m({P_{\ell}},k-1)+m({P_{\ell-2}},k-2)\\
&=&m({{C_{\ell}}},k-1)\\
&\leq &
b_{2k-2}({\overrightarrow{C_{\ell}}}^+),\;\;\;\;\;\;\; \mbox( by \
Corollary  \ \ \ref{bre})
\end{eqnarray*}
where $s+t=l-2$.

Suppose ${{\G}^+}\prec {\overrightarrow{P_{n'}^{\ell}}}^+$ for all
$n'<n$.  Then there is a pendant edge, say $uv$, with pendant
vertex $v$. By Lemma \ref{cut},

\begin{eqnarray*}
b_{2k}({{\G}^+})&=&b_{2k}({{\G}^+}-v)+
b_{2k-2}({{\G}^+}-v-u);\\
b_{2k}({\overrightarrow{P_n^{\ell}}}^+)&=&b_{2k}({\overrightarrow{P_{n-1}^{\ell }}}^+)+b_{2k-2}({\overrightarrow{ P_{n-2}^{\ell}}}^+).
\end{eqnarray*}

By induction hypothesis, it  suffices to prove that
$b_{2k-2}({{\G}^+}-v-u)\leq b_{2k-2}({\overrightarrow{
P_{n-2}^{\ell}}}^+)$. If ${{\G}^+}-v-u$ contains a cycle, then by
induction hypothesis the inequality holds. Suppose that
${{\G}^+}-v-u$ is a forest, then  by Lemmas \ref{match} and
\ref{coeff},
\begin{eqnarray*}
b_{2k-2}({{\G}^+}-v-u )&=& m({{G}}-v-u,k-1)\\
&\leq & m({P_{n-2}},k-1)\leq
m({{P_{n-2}^{\ell}}},k-1)\\ &\leq&
b_{2k-2}({{ P_{n-2}^{\ell}}}). \mbox{ \hspace*{1cm} $\Box$ } \end{eqnarray*}

\begin{Lem} For $3\leq \ell \leq n$, $\ell\neq 4$, ${\overrightarrow{P_n^{\ell}}}\prec {\overrightarrow{P_n^{4}}}^+$. \label{pnl}\end{Lem}

\begin{Proof} We proceed the proof by induction on $n$. For $n=\ell=5$, $b_4({\overrightarrow{P_5^{5}}}^+)=
b_4({\overrightarrow{C_5}}^+)=5$,
$b_4({\overrightarrow{P_5^4}}^+)=b_4({\overrightarrow{S_5^4}}^+)=6$,
thus $b_4({\overrightarrow{P_5^{5}}}^+)<
b_4({\overrightarrow{P_5^4}}^+)$. Suppose that
${\overrightarrow{P_{n'}^{\ell}}}^+\prec
{\overrightarrow{P_{n'}^{4}}}^+$ for all $n'<n$,  $5\leq \ell$.
By Lemma \ref{char}, we have
\begin{eqnarray*}b_{2k}({\overrightarrow{P_n^{\ell}}}^+)&=& b_{2k}({\overrightarrow{P_{n-1}^{\ell}}}^+)+b_{2k-2}({\overrightarrow{ P_{n-2}^{\ell}}}^+), \\
b_{2k}({\overrightarrow{P_n^{4}}}^+)&=& b_{2k}({\overrightarrow{P_{n-1}^{4}}}^+)+b_{2k-2}({\overrightarrow{ P_{n-2}^{4}}}^+).
\end{eqnarray*}

By induction hypothesis, the result follows.

Similarly, we can prove that for $\ell=3$ the inequality  holds too.
\end{Proof}

By Lemmas \ref{gless} and  \ref{pnl}, we have

\begin{Thm} Among all  orientations of unicyclic graph,  ${\cal E}_s({\overrightarrow{P_n^{4}}}^+)$  is the unique oriented graph (under switching-equivalent)  with maximal skew energy.  \end{Thm}

{\bf Acknowledgments.}
This project  is supported
by National Natural Science Foundation of China


\begin{thebibliography}{01}

\footnotesize


\bibitem{ca} C. Adiga, R. Balakrishnan and Wasin So, The skew energy of a digraph, {\em Linear Algebra
and its Application} 432 (2010) 1825-1835.

\bibitem{dmc} D. M. Cvetkovic, M. Doob, H. Sachs, {\em Spectra of Graphs,} Academic Press, New York, 1979.

\bibitem{md} M. Dehmer, F. Emmert-Streib, Analysis of complex networks: from biology to
linguistics. WILEY-VCH Verlag GmbH £¦ Co. hGaA, Weinheim, 2009.

\bibitem{gutman} I. Gutman, The energy of a graph, {\em Ber. Math. Stat. Sekt. Forschungszentrum Graz} 103 (1978) 1--22.

\bibitem{gu} I. Gutman, Y. Hou,  Bipartite unicyclic
graphs with greatest energy, {\em MATCH Commun. Math. Comput. Chem.} 43
(2001) 17--28.

\bibitem{bh} B. Huo, X. Li, Complete solution to a conjecture on
the maximal energy of unicyclic graphs. {\em European journal of
combinatorics.} 32 (2011) 662--673.

\bibitem{hh2}  H. Hua,  M. Wang, Unicyclic graphs
with given number of pendent vertices and minimal energy, {\em Lin.
Algebra Appl.} 426 (2007) 478--489.

\bibitem{hou1} Y. Hou, Unicyclic graphs with minimal energy, {\em J. Math. Chem.}
29 (2001) 163--168.
\bibitem{hou2} Y. Hou, I. Gutman, C. W. Woo,  Unicyclic
graphs with maximal energy, {\em Lin. Algebra Appl.} 356 (2002) 27--36.
\bibitem{hou} Y. Hou, T. Lei, Characteristic polynomials of skew-adjacency matrices of oriented
graphs, The Electronic Journal of Combinatorics 18(2011) \# P156.



\bibitem{lx}  X. Li, J. Zhang,  B. Zhou, On unicyclic
conjugated molecules with minimal energies, {\em J. Math. Chem.} 42
(2007) 729--740.
\bibitem{vn} V. Nikiforov, The energy of graphs and matrices, {\em J. Math. Anal. Appl.} 320,(2007) 1472--1475.

    \bibitem{pr} I. Pena, J. Rada, Energy of digraphs, Linear and Multilinear Algebra 56(2008) 565-579.


\bibitem{whw} W. H. Wang,  A. Chang,  L. Z. Zhang,  D. Q. Lu, Unicyclic H$\ddot{u}$ckel
molecular graphs with minimal energy, {\em J. Math. Chem.} 39 (2006),
231--241.

\bibitem{whw2} W. H. Wang,   A. Chang,  D. Q. Lu, Unicyclic
graphs possessing Kekul¨¦ structures with minimal energy, {\em J. Math.
Chem.} 42 (2007) 311--320.




\end{thebibliography}
\end{document}